\documentclass{ifacconf}

\usepackage{graphicx,float}      
\usepackage{natbib}        
\usepackage{tikz}
\usepackage{amsfonts,latexsym,amssymb} 
\usepackage{mathrsfs,amsmath}

\newcommand{\vecA}{\mathbf{A}}
\newcommand{\vecB}{\mathbf{B}}
\newcommand{\vecC}{\mathbf{C}}

\newcommand{\vecF}{\mathbf{F}}
\newcommand{\vecG}{\mathbf{G}}
\newcommand{\vecH}{\mathbf{H}}
\newcommand{\vecI}{\mathbf{I}}

\newcommand{\vecM}{\mathbf{M}}

\newcommand{\vecT}{\mathbf{T}}

\newcommand{\vecX}{\mathbf{X}}

\newcommand{\vecb}{\mathbf{b}}
\newcommand{\vecc}{\mathbf{c}}
\newcommand{\vece}{\mathbf{e}}
\newcommand{\vecf}{\mathbf{f}}
\newcommand{\vecg}{\mathbf{g}}

\newcommand{\vecu}{\mathbf{u}}

\newcommand{\vecx}{\mathbf{x}}

 \newtheorem{proposition}{Proposition}
     \newtheorem{remark}{Remark}[section]
      
        \newtheorem{lemma}{Lemma}
\begin{document}
\begin{frontmatter}

\title{Purcell magneto-elastic swimmer controlled by an external magnetic field}


\author[First]{F. Alouges} 
\author[Second]{A. DeSimone} 
\author[Third]{L. Giraldi}
\author[Forth]{M. Zoppello}

\address[First]{CMAP UMR 7641 \'Ecole Polytechnique CNRS, Route de Saclay, 91128 Palaiseau Cedex France (e-mail: francois.alouges@polytechnique.edu)}
\address[Second]{Scuola Internazionale Superiore di Studi Avanzati via Bonomea 265 I-34136 Trieste Italy (e-mail: desimone@sissa.it)}
\address[Third]{Universit\'e C\^ote d'Azur, CNRS, LJAD, INRIA Sophia Antipolis Méditerranée Team/équipe McTAO B.P. 93 -- 06902 Sophia Antipolis cedex -- France (e-mail: laetitia.giraldi@inria.fr) }
\address[Forth]{Universit\'a degli studi di Trento, via Sommarive 14 - 38123 Povo (TN) Italy (e-mail: marta.zoppello@unitn.it)}

\begin{abstract}                
This paper focuses on the mechanism of propulsion of a Purcell swimmer whose segments are magnetized and react to an external magnetic field applied into the fluid. 
By an asymptotic analysis, we prove that it is possible to steer the swimmer along a chosen direction when the control functions are prescribed as an oscillating field. 
Moreover, we discuss what are the main obstructions to overcome in order to get classical controllability result for this system.
\end{abstract}


\end{frontmatter}

\section{Introduction}
In the last decade, Micro-motility has become a subject of growing interest, both for the biological understanding of micro-organisms and technological applications. In the latter direction, the topic addresses several challenges as for instance the conception of artificial self-propelled and/or easily controllable microscopic robots. Such kind of devices could revolutionize the biomedical applications \cite{PeyerZhang05} as for instance it could be useful to minimize invasive microsurgical operations \cite{NelsonKaliakatsos10}.  
\\

One of the few possibilities recently studied in the literature, 
is to consider an artificial swimmer that possesses a magnetic flexible tail, and use an external magnetic field to act on this flagellum \cite{Bibette,GaoSattayasamitsathit10,GaoKagan12,PakGao11}. 
This particular design is influenced by the locomotion of spermatozoa which achieves their propulsion by propagating travelling wave along their flagellum. \\

On the other hand, there exists now a quite wide literature that makes a connection between problem of swimming at the micro-scale, and the mathematical control theory. Starting from the pioneering work of Shapere and Wilczek \cite{Shapere89}, and Montgomery \cite{Montgomery02}, the dynamics of self-propelled microscopic artificial swimmers has been considered for instance in \cite{AlougesDeSimone13, AlougesDeSimone08,  AlougesGiraldi12, Gerard-VaretGiraldi13}) where the rate of shape changes of the swimmer is considered as a natural control. 

The aim of the present paper is therefore to study thoroughly the mathematics of a magneto-elastic artificial swimmer, as the ones proposed in \cite{AlougesDeSimone14}, when controlled with an external oscillating magnetic field. 
As in \cite{Or14}, we provide an estimate of the displacement of the swimmer with respect to small amplitude of the fields, generalizing the approach of E. Gutman and Y. Or, to the case of a $3$-link swimmer.\footnote{The 3-link swimmer was already proposed in \cite{Purcell77}. It is one of the more famous systems, which was further considered in the literature \cite{BeckerKoehler03,GiraldiMartinon13,Or,TamHosoi07} for instance.}  

In our model the links are supposed to be uniformly magnetized and linked together with rotational springs. When an external magnetic field is applied its shape is expected to experience a deformation which hopefully leads to a global displacement. It was shown in \cite{AlougesDeSimone14} that the dynamics of the a magneto-elastic swimmer is governed by a system of ODEs which is an affine control system with respect to the magnetic field. 
We adopt the same approach of \cite{AlougesDeSimone14} to get the equations of motion. At the micro-scale, the flow is characterized by small Reynolds number. Thus, we assume that the surrounding fluid is governed by Stokes equations which implies that hydrodynamic forces and torques are linear with respect to the  velocity distribution on the boundary (rates of deformation and displacement). In that case, the Resistive Force Theory (RFT) (see \cite{GrayHancock55}) provides a simple and concise way to compute a local approximation of hydrodynamic forces involved in our system.
The magnetic  behavior of the segments is modeled by assuming that their magnetization is always parallel to the segment with fixed magnitude and stray fields, especially magnetic interactions between different segments, are neglected. Only the magnetic torque induced by the external magnetic field on each segment is considered.  
\\

After briefly recalling the equation of motion in Section \ref{2}, we provide, in Section \ref{3} an estimate of the displacement of the swimmer for an oscillating magnetic field of small amplitude, which steers it along one direction. The last section \ref{4} of the paper is devoted to discuss the main obstruction to get classical controllability results.

\section{Modeling}
\label{2}
The model of the magneto-elastic $N$-link swimmer was already introduced in \cite{AlougesDeSimone14}. The two link one was studied in \cite{Or14, GiraldiPomet15}.
Here, we recall briefly the equation of motion focusing on the case $N=3$.

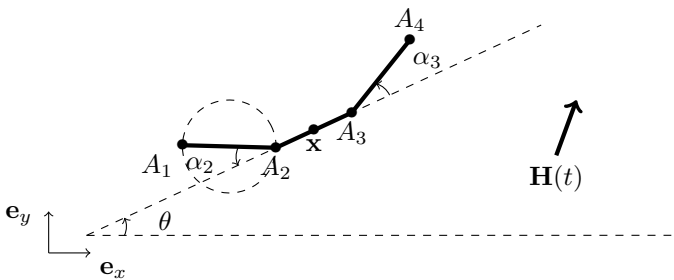
\begin{figure}[H]
 \begin{center}
\begin{tikzpicture}[scale=0.55,rotate=25]
\draw[dashed] (-6,0)--(6,0);
\draw[ultra thick] (-1,0) -- (1,0);
\draw[ultra thick] (-1,0) -- ( -3,1);
\draw[ultra thick] (1,0) -- ( 3,1);
\draw [dashed] (-6,0) -- (6.68,-5.91);
\draw [->](-7,0) -- (-7 +0.906,-0.422);
\draw [->] (-7,0) -- (-7 +0.422, 0.906);
\draw[ultra thick] [->] (5,-3) -- (6,-2);
 \draw[dashed] (-2,0.5) circle (1.12cm);
 \draw (3,1) node {$\bullet$};
   \draw (0,0) node {$\bullet$};
  \draw (-1,0) node{$\bullet$};
  \draw (1,0) node{$\bullet$};
  \draw (-3,1) node{$\bullet$};
   \draw [->](2,0) arc (0:50:0.5);
   \draw [<-] (-2,0) arc (0:-50:-0.5);
   \draw [<-] (-5,0) arc (0:-50:0.5);
\draw (-3,1) node[below left] {$\tiny{A_1}$};
 \draw (-1,0) node[below] {$\tiny{A_2}$};
  \draw (1,0) node[below] {$\tiny{A_3}$};
  \draw (3,1) node[above] {$\tiny{A_4}$};
   \draw (0,0) node[below] {$\tiny{\vecx}$};
     \draw (-3,-0.) node[above] {$\alpha_2$};
      \draw (3,-0.) node[above] {$\alpha_3$};
      \draw (-4.3,0.1) node[below right] {$\theta$};
       \draw (-7 +0.906,-0.422) node[below right] {$\tiny{\vece_x}$};
        \draw (-7 +0.422, 1.3)node[below left] {$\tiny{\vece_y}$};
        \draw (5,-3) node[below] {$\tiny{\vecH(t)}$};
     \end{tikzpicture}
     \caption{The magneto-elastic Purcell swimmer of shape $(\alpha_2,\alpha_3)$ at the position $\vecx$ in the plane subject to an external magnetic fields $\vecH$.}
 \label{3_links_swimmer}
     \end{center}
     \end{figure}

We consider a magneto-elastic 3-link Purcell swimmer moving in a plane. 
This two-dimensional setting is suitable for the study of slender, essentially one-dimensional swimmers exploring planar trajectories as explained in \cite{AlougesDeSimone13, AlougesDeSimone14}.
The swimmer consists of $3$ rigid segments, each of length $L$ with articulated joints at their ends (see Fig. \ref{3_links_swimmer}), moving in the horizontal $2$d-plane of the lab-frame.
Because of the symmetric geometry of the swimmer, we use slightly different notation and variables than in \cite{AlougesDeSimone14}. Indeed we call $\vecx=(x,y)$ the coordinates of the central point of the second segment, $\theta$ the angle that it forms with the $x$-axis, $\alpha_2$ the relative angle between the first and second segments and finally $\alpha_3$ the relative between the third and the second segments (see Fig \ref{3_links_swimmer}). Therefore the position and the orientation of the swimmer are characterized by the triplet $(x,y,\theta)$, while its shape is given by $(\alpha_2,\alpha_3)$.

As in \cite{AlougesDeSimone14} the three segments are uniformly magnetized and linked together with torsional springs, with elastic constant $K$, that tend to align the segments one with another. Those produce torques when the segments are not fully aligned. In what follows we assume that $\vecH(t) :=\begin{pmatrix}H_x(t)\\H_y(t)\\0 \end{pmatrix}$ is horizontal in such a way that the motion holds in the horizontal plane.

\subsection{Equations of motion}
\label{Dynamic_of_the_swimmer}
As it was noticed in \cite{AlougesDeSimone14} the equations which govern the dynamics of the swimmer form a system of ODEs, which is affine with respect to the magnetic field $\vecH(t)$. 
The hydrodynamic forces acting on the $i$-th link, are approximated using RTF, with parallel (resp. perpendicular) drag coefficients $\xi_i$ (resp. $\eta_i$). In our particular case this system describes the evolution of the position and the shape variables and thus consists of five equations. Those are obtained by writing the balance of forces for the whole swimmer and the balance of torques for the subsystems consisting of the three, two and one rightmost segments. We call those subsystems $S_1,S_2,S_3$ respectively ($S_1$ is therefore the whole system). The motion of the swimmer holds in the horizontal plane since only horizontal forces and vertical torques apply. The final system reads:

\begin{equation}
\label{Newton_laws}
\left\{
\begin{array}{lllll}
\vecF_h = 0 \,,\\
\vece_z \cdot \left(\vecT^{A_1}_h +\vecT^{A_1}_e+ \vecT^{A_1}_m\right) =0\,,\\
\vece_z \cdot \left(\vecT^{A_2}_h + \vecT^{A_2}_e + \vecT^{A_2}_m\right) = 0\,,\\
\vece_z \cdot \left(\vecT^{A_3}_h+ \vecT^{A_3}_e + \vecT^{A_3}_m\right)=0\,.
\end{array}
\right.
\end{equation}

Here, $\vecF_h$ denotes the total hydrodynamic force acting on the swimmer, $\vecT^{A_i}_h$ (resp. $\vecT^{A_i}_e$ and $\vecT^{A_i}_m$) is the hydrodynamic (resp. elastic  and magnetic) torque with respect to $A_i$ (see Fig. \ref{3_links_swimmer}), acting on the subsystem $S_i$.
Following the construction made in \cite{AlougesDeSimone14} system \eqref{Newton_laws} becomes
\begin{equation}
\begin{aligned}
\label{eq:all_system}
&\vecM_h \left(\theta,\alpha_2,\alpha_3\right)\begin{pmatrix}
\dot{\vecx}\\
\dot{\theta}\\
\dot{\alpha}_2\\
\dot{\alpha}_3
\end{pmatrix}
= -K \begin{pmatrix} 0 \\ 0 \\ 0\\ \alpha_2\\\alpha_3\end{pmatrix}\\&\\&
-\vecM^x_m \left(\theta,\alpha_2,\alpha_3\right)H_x(t) -\vecM^y_m \left(\theta,\alpha_2,\alpha_3\right)H_y(t)\,,
\end{aligned}
\end{equation}
with $\vecM_h$ $5\times 5$ matrix, $\vecM^x_m $ and $\vecM^y_m$ vectors in $\mathbb{R}^5$ all depending on $(\theta,\alpha_2,\alpha_3)$. All these matrices can be computed explicitly following the approach given in \cite{AlougesDeSimone14}.\\
%
%
Finally system \eqref{eq:all_system} can be rewritten as 
\begin{equation}
\begin{aligned}
\label{eq:full_system}
\begin{pmatrix}
\dot{\vecx}\\
\dot{\theta}\\
\dot{\alpha}_2\\
\dot{\alpha}_3
\end{pmatrix}
=&\vecf_0(\theta,\alpha_2,\alpha_3)+\vecf_x(\theta,\alpha_2,\alpha_3)H_x(t)\\&+\vecf_y(\theta,\alpha_2,\alpha_3)H_y(t)
\end{aligned}
\end{equation}
where,
 \begin{equation}
 \label{eq:def_of_f}
 \begin{array}{ll}
 \vecf_0= &-\vecM_h^{-1} K \begin{pmatrix} 0 \\ 0 \\ 0\\ \alpha_2\\\alpha_3\end{pmatrix}\\
\vecf_j=&-\vecM_h^{-1}\vecM^j_m\,,\quad\text{ $j=x,y$. } \\
  \end{array}
  \end{equation}
 
Notice however that the dynamics of $(\theta,\alpha_2,\alpha_3)$ is independent of $\vecx$ and can be decoupled. Indeed by block decomposing the matrix $\vecM_h$ as
\begin{equation}
\label{eq:matrix_block}
\vecM_h=\begin{pmatrix}
\vecA_h&\vecB_h\\
\vecB^T_h&\vecC_h\end{pmatrix}
\end{equation}
($\vecA_h$, $\vecB_h$ and $\vecC_h$ being respectively $2\times 2$, $2\times 3$ and $3\times 3$ matrices), and considering the first two rows of the system \eqref{eq:all_system}, we can solve for $\dot{\vecx}$ as
\begin{equation}
\small
\label{eq:position}
\begin{array}{ll}
\dot{\vecx}&=-\vecA_h^{-1}\left(\theta,\alpha_2,\alpha_3\right)\vecB_h\left(\theta,\alpha_2,\alpha_3\right)\begin{pmatrix}
\dot{\theta}\\
\dot{\alpha}_2\\
\dot{\alpha}_3
\end{pmatrix}\\
&= \vecG_1\left(\theta,\alpha_2,\alpha_3\right) \dot{\theta} + \vecG_2\left(\theta,\alpha_2,\alpha_3\right) \dot{\alpha_2} +\vecG_3\left(\theta,\alpha_2,\alpha_3\right) \dot{\alpha_3}\,.
\end{array}
\end{equation}

Moreover, the subsystem associated with the shape and the orientation of the swimmer becomes
\begin{equation}
\begin{aligned}
(-\vecB_h^T&\vecA^{-1}_h\vecB_h +\vecC_h)\begin{pmatrix}
\dot{\theta}\\
\dot{\alpha}_2\\
\dot{\alpha}_3
\end{pmatrix} =\\&-K \begin{pmatrix} 0\\ \alpha_2\\\alpha_3\end{pmatrix}
-
\tilde{\vecM}^x_m H_x(t) -\tilde{\vecM}^y_m H_y(t)\,,
\end{aligned}
\end{equation}
that we rewrite inverting the left hand side matrix as
\begin{equation}
\begin{aligned}
\label{eq:angles}
\begin{pmatrix}
\dot{\theta}\\
\dot{\alpha}_2\\
\dot{\alpha}_3
\end{pmatrix} &= \vecg_0(\theta,\alpha_2,\alpha_3) +\\&+\vecg_x(\theta,\alpha_2,\alpha_3) H_x(t) + \vecg_y(\theta,\alpha_2,\alpha_N) H_y(t) \,.
\end{aligned}
\end{equation}

The whole dynamical system \eqref{eq:full_system} is an affine control system with drift where the controls are the two components of the magnetic field. The explicit expression of the dynamics are formally computed by using a symbolic computation software as Mathematica.

\section{Steering along one direction with small sinusoidal magnetic fields}
\label{3}
As it was noticed by experiments in \cite{PoperDreyfus06}, in the rest we show,
that a swimmer with an initial symmetric shape, remains symmetric. 
Moreover by using an asymptotic expansion we provide an estimate of the swimmer displacement for a prescribed small sinusoidal magnetic field.
\subsection{A symmetry obstruction}
\label{subsec:symmetric_obstruction}

In this part, we assume that the drag coefficients of the $3$ links are identical and that they have all a uniform magnetization,  
$$\xi_i=\xi \quad \text{and} \quad \eta_i=\eta\,,\quad M_i=M\,\, \forall i\in \left\{1,2,3\right\}\,.$$ 
These assumptions are suitable to describe the behavior of a magnetic filament.\\
We consider the set of symmetric shapes  $$\mathcal{S}:=\left\{(\alpha_2,\alpha_3)\in[0,2\pi]^2, \alpha_2=\alpha_3 \right\},$$ (see Fig. \ref{Fig_symm}). 
As it was experimentally observed in \cite{PoperDreyfus06}, the following Proposition \ref{prop:symmetric} proves that for a swimmer which possesses an initial shape belonging to $\mathcal{S}$, its shape remains symmetric, regardless of the applied magnetic field.

\begin{proposition}
\label{prop:symmetric}
Let $T>0$, if at the initial time $\alpha_2(0)=\alpha_3(0)$ then for any magnetic field  $t\mapsto (H_x(t),H_y(t))$, applied to the system, the shape of the swimmer remains symmetric i.e., 
$$\alpha_2(t)=\alpha_3(t)\,, \forall t\in[0,T]\,.$$
\end{proposition}

\begin{pf}
This result is based on a symmetry argument. The dynamics must be invariant by the rotation $\mathcal{R}$ of angle $\pi$  about the center $O$ of the lab-frame of the whole system (swimmer and magnetic field). Notice that this rotation changes $\vecx$ to $-\vecx$, leaves $\theta$ invariant, interchanges $\alpha_2$ and $\alpha_3$, and reverses the magnetic field $\vecH$ and the magnetization along the swimmer. Therefore, if the function $t\mapsto (\vecx(t),\theta(t),\alpha_2(t),\alpha_3(t))$ is solution of the system  \eqref{eq:full_system} magnetized in one direction (say $+$), for an external magnetic field $\vecH(t)$, then the trajectory $(-\vecx(t),\theta(t),\alpha_3(t),\alpha_2(t))$ is the solution corresponding to the magnetic field $-\vecH(t)$ of the system magnetized in the opposite direction (say $-$). We summarize this by saying that
\begin{equation}
\footnotesize
(\vecx(t),\theta(t),\alpha_2(t),\alpha_3(t),\vecH(t),+) \xrightarrow{\mathcal{R}} (-\vecx(t),\theta(t),\alpha_3(t),\alpha_2(t),-\vecH(t),-)
\end{equation}
where the last component corresponds to the direction of the magnetization along the swimmer.
 
Similarly, we consider a second transformation $\mathcal{T}$ which reverses only the magnetization and the external magnetic field. We remark that the equations of motion \eqref{eq:full_system} depend only on the product $M\vece_k\wedge\vecH(t)$, if  the function $t\mapsto (\vecx(t),\theta(t),\alpha_2(t),\alpha_3(t))$ is solution of \eqref{eq:full_system} with a prescribed magnetic field and magnetization, then it remains a solution with opposite magnetic field and magnetization, so that
\begin{equation}
\footnotesize
(\vecx(t),\theta(t),\alpha_2(t),\alpha_3(t),\vecH(t),+) \xrightarrow{\mathcal{T}} (\vecx(t),\theta(t),\alpha_2(t),\alpha_3(t),-\vecH(t),-)\,.
\end{equation}

The geometric transformations $\mathcal{R}$ and $\mathcal{T}$ are sketched in Fig. \ref{Fig_symm}. Checking more formally those symmetry properties of the systems can of course be done on \eqref{eq:all_system} but is left to the reader.

Now, composing $\mathcal{R}$ and $\mathcal{T}$ we have, using the notation above
\begin{equation}
\footnotesize
(\vecx(t),\theta(t),\alpha_2(t),\alpha_3(t),\vecH(t),+) \xrightarrow{\mathcal{T}\circ \mathcal{R}} (-\vecx(t),\theta(t),\alpha_2(t),\alpha_3(t),\vecH(t),+)\,,
\end{equation}
which means that by uniqueness of the solution of \eqref{eq:full_system}, a swimmer starting at position $\vecx(0)=0$ with a symmetric shape ($\alpha_2(0)=\alpha_3(0)$),
under whatever driving magnetic field $\vecH(t)$ verifies
$$
\vecx(t) = -\vecx(t)\,, \mbox{ and }\alpha_2(t) = \alpha_3(t)\,.
$$
It hence experiences no displacement and stays symmetric. 
%




\begin{figure}[H]
\begin{center}
\includegraphics[width=8.4cm]{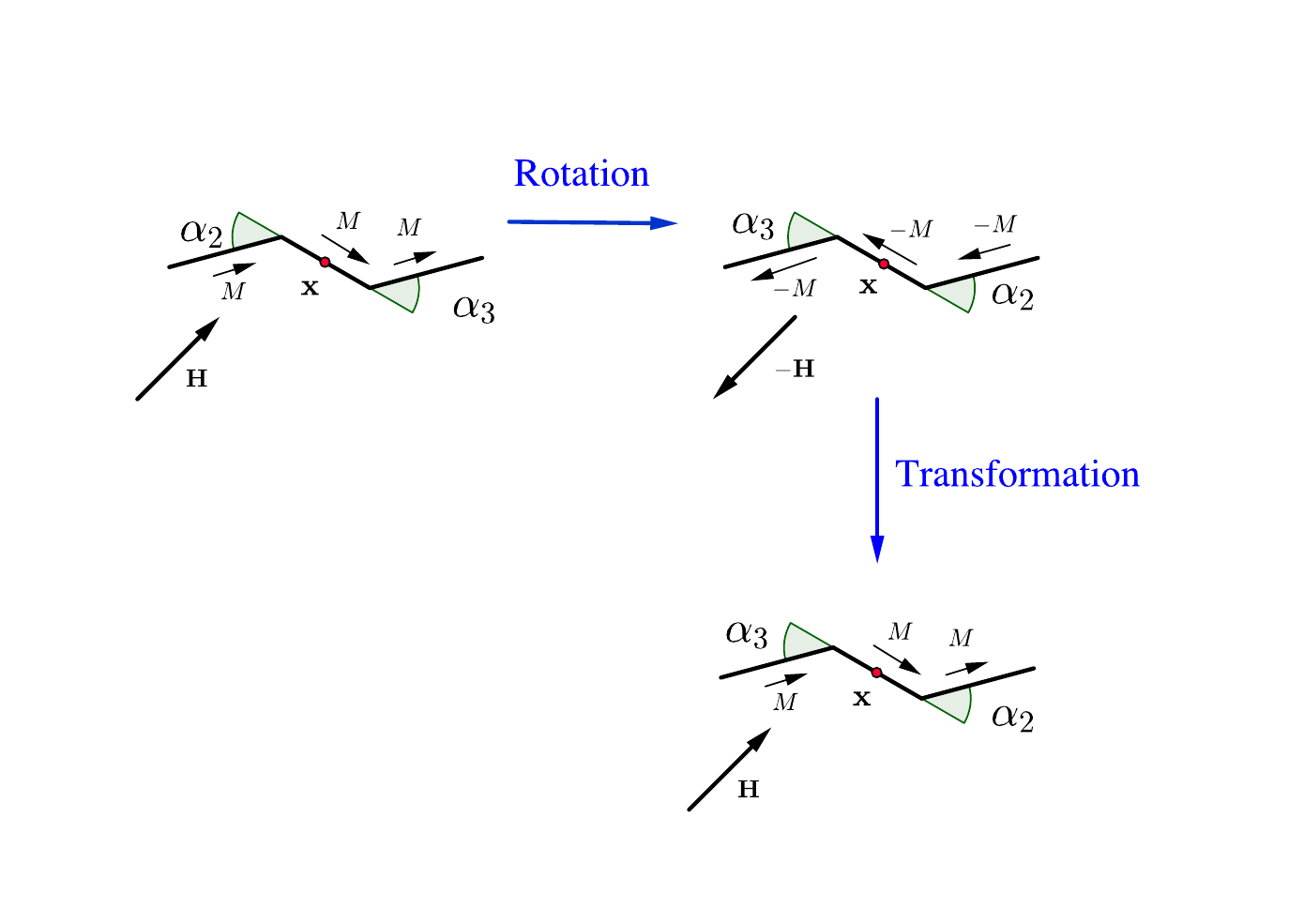}
\caption{The two transformations $\mathcal{R}$ and $\mathcal{T}$ of the system.\label{Fig_symm}}
\end{center}
\end{figure}
\end{pf}
\begin{remark}
Similar argument holds for a swimmer composed by an odd number of links.
\end{remark}
\subsection{Small oscillating magnetic field}
In this subsection we focus on a swimmer satisfying 
\begin{equation}
\label{drag_coeff_ass}
\begin{aligned}
&\eta_2=\eta_3=\eta&&\xi_2=\xi_3=\xi\\
& \text{with   } \eta_1\neq \eta&&\xi_1\neq\xi
\end{aligned}
\end{equation}
These latter assumptions allow to overcome the previous symmetry obstruction. They are suitable in the case of a swimmer with an head.\\
Starting from a swimmer
with a horizontal shape $(\theta=\alpha_2=\alpha_3=0)$, we can use
the horizontal component of the external magnetic field as a ``stabilizer''
whereas the oscillating vertical component produces the shape deformation 
and the motion. In order to understand further what happens when such
a field is applied we make the following perturbation analysis. 
We assume that
\begin{equation}
\label{hyp:particular_field}
(H_x(t),H_y(t))=(1,\epsilon\sin(\omega t))\,,
\end{equation}
and compute the asymptotic expansion of the swimmer displacement with respect to small $\epsilon$ after a period $\frac{2\pi}{w}$. 

Linearizing the system of equations  \eqref{eq:angles} for small angles $(\theta,\alpha_2,\alpha_3)$, to first order in $\epsilon$, i.e.,
$$
\begin{pmatrix}
\theta\\
\alpha_2\\
\alpha_3
\end{pmatrix} = \epsilon \begin{pmatrix} \tilde{\theta}\\\tilde{\alpha}_2\\\tilde{\alpha}_3\end{pmatrix}+ o(\epsilon)
$$
we get that the triplet $\begin{pmatrix} \tilde{\theta}\\\tilde{\alpha}_2\\\tilde{\alpha}_3\end{pmatrix}$ satisfies the equation, 

\begin{equation}
\label{eq:angles_small}
\begin{pmatrix} \dot{\tilde{\theta}}\\\dot{\tilde{\alpha}}_2\\\dot{\tilde{\alpha}}_3\end{pmatrix} = \vecA \begin{pmatrix} \tilde{\theta}\\\tilde{\alpha}_2\\\tilde{\alpha}_3\end{pmatrix}+b \sin(\omega t)
\end{equation}
with
\begin{equation}
\vecA= \nabla \left(\vecg_0 + \vecg_x\right) \left(0,0,0\right)\,,\,\, b=\vecg_y(0,0,0)\,.
\end{equation}
Here, $\vecA$ is the $3\times3$ matrix which depends on the drag coefficients, $(\eta_1,\eta)$ and $(\xi_1,\xi)$, on the magnetization $M$ and on the elastic constant $K$. We find that its explicit expression is given by

%
%
\begin{equation*}
\small
\vecA=\delta
\begin{pmatrix}
a_{1\,1}&a_{1\,2}&a_{1\,3}\\
a_{2\,1}&a_{2\,2}&a_{2\,3}\\
a_{3\,1}&a_{3\,2}&a_{3\,3}
\end{pmatrix}
\end{equation*}

where
\begin{align*}
 &a_{1\,1}=M\eta_1(5\eta+\eta_1) \\
 &a_{1\,2}= (19 K+9M) \eta\eta_1 +2 K \eta_1^2\\
 &a_{1\,3}= 2 (8 K+3M) \eta\eta_1 +(5 K+3M) \eta_1^2\\
 &a_{2\,1}= -M\left( 4 \eta ^2 +13 \eta_1 \eta +\eta_1^2\right)\\
 &a_{2\,2}= -4 (K+M) \eta^2-(42 K+23M) \eta_1 \eta -2 K \eta_1^2\\
 &a_{2\,3}= - (28 K+9M) \eta \eta_1 -(5 K+3M) -\eta_1^2\\
 & a_{3\,1}= -6M (2 \eta\eta_1 +\eta_1^2)\\
 &a_{3\,2}= - 4 (7 K+3M) \eta \eta_1-5 K \eta_1^2 \\
 &a_{3\,3}= - 16 (2 K+M) \eta \eta_1-(16 K+11M) \eta_1^2 
\end{align*}

and
$\delta=\frac{6}{L^3 \eta \eta_1 (8 \eta +7 \eta_1)}$.
\newline


The solution of the system \eqref{eq:angles_small} is thus given by 
\begin{equation}
\small
\begin{aligned}
\label{eq:sol_edo_angles}
\begin{pmatrix} \tilde{\theta}(t)\\\tilde{\alpha}_2(t)\\\tilde{\alpha}_3(t)\end{pmatrix} &=\frac{1}{2i}\left( \vecA^+(\omega) \exp{(i\omega t)} b - \vecA^-(\omega) \exp{(-i\omega t)} b\right) \\&+ \vecc(\omega)  \exp{(\vecA t)}
\end{aligned}
\end{equation}
where $ \vecA^{\pm}(\omega) := (-\vecA \pm i \omega\vecI )^{-1}$ and $\vecc(\omega) = \vecA^{-}(\omega)b - \vecA^{+}(\omega)b$.
\\

The first part of the solution corresponds to a periodic solution, while the last is an
exponentially decaying perturbation, as we shall see now. Indeed, by applying Routh-Hurwitz criterion on the characteristic polynomial of $\vecA$, we prove that the real part of its eigenvalues are all negative. This provides the stability of the asymptotic periodic solution 
\\

Let us recall the Routh-Hurwitz criterion.

\begin{lemma}[Routh-Hurwitz criterion, \cite{Gantmacher59}]
\label{lemma:Hurwitz}
Let $P(s) = a_3s^3 + a_2s^2 + a_1s + a_0 = 0$ a polynomial of third degree. If the two following conditions are satisfied
\begin{enumerate}
\item All the coefficients have the same sign,  
\item $a_2a_1 > a_3a_0$,  
\end{enumerate}
then the real part of the roots of $P$ are strictly negative. 
\end{lemma}

Here the characteristic polynomial of $\vecA$ reads
  $$p_{\vecA}(\lambda) = \mbox{det}(\vecA-\lambda \textbf{Id}) = a_3\lambda^3 + a_2 \lambda^2 + a_1 \lambda + a_0\,,$$
where
\begin{align*}
\small
a_0&=\mbox{det}\vecA=-\frac{432 M\left(3 K^2+4 KM+M^2\right) (2 \eta +\eta_1)}{d}<0\,,\\
a_1&=-\frac{36}{d}\bigg( (M^2\left(10\eta^2+28\eta\eta_1+\eta_1^2\right)\\
&\phantom{= -\frac{36}{d} }+K^2 \left(16\eta^2+64\eta\eta_1+\eta_1^2\right)\\
&\phantom{= -\frac{36}{d} }+KM\left(31\eta^2+98\eta\eta_1+3\eta_1^2\right) )\bigg)<0\,,\\
a_2&=\mbox{tr}\vecA=\frac{-12}{d}\bigg(M\left(2\eta^2+17\eta\eta_1+5\eta_1^2\right)\\&\phantom{\mbox{tr}\vecA=\frac{-12}{d}}+K\left(2\eta^2+37\eta\eta_1+9\eta_1^2\right)\bigg)<0\,,\\
a_3 &=-1<0\\
d&=L^9 \eta ^2\ \left(8 \eta  \eta_1+7 \eta_1^2\right)
\end{align*}

so that condition 1. of Lemma \ref{lemma:Hurwitz} is satisfied. Moreover,
\begin{equation}
\small
\begin{aligned}
a_3a_0-&a_2a_1=\frac{-1}{L^9 \eta ^3 \eta_1^2 (8 \eta +7 \eta_1)^2}\Bigg(432 \bigl(2(K+M)(16 K^3\\&+31 KM+10M^2) \eta ^4+\\
&+5\left(144 K^3+339 K^2M+217KM^2+42M^3\right) \eta ^3 \eta_1\\
&+\left(2514 K^3+5015 K^2M+2867 KM^2+506M^3\right) \eta ^2 \eta_1^2+\\
&+\left(613 K^3+1309 K^2M+802 KM^2+150M^3\right) \eta  \eta_1^3\\
&+(9 K+5M)\left(K^2+3KM+M^2\right)\bigr) \eta_1^4\Bigg)<0\,.\\
\end{aligned}
\end{equation}
Therefore, we deduce that the steady state of the equation \eqref{eq:angles_small} is stable, i.e.,
\begin{equation}
\label{eq:steady_state}
\exp{(\vecA t)} \to 0\,, \quad \mbox{ as } \quad t \to \infty\,,
\end{equation}
\newline
and the solution of \eqref{eq:sol_edo_angles} exponentially converges to the periodic solution
\begin{equation}
\label{eq:sol_edo_angles_infty}
\begin{pmatrix} \tilde{\theta}^\infty(t)\\\tilde{\alpha}^\infty_2(t)\\\tilde{\alpha}^\infty_3(t)\end{pmatrix} =\frac{1}{2i}\left( \vecA^+(\omega) \exp{(i\omega t)} b - \vecA^-(\omega) \exp{(-i\omega t)} b\right)\,,
\end{equation}
and in particular $\theta \sim \epsilon \tilde{\theta}^\infty$ oscillates around 0 indicating that the swimmer stays nearly horizontal, stabilized by the horizontal component of the magnetic field. Similarly, the fact that the shape variables $(\alpha_2,\alpha_3)$ are periodic (and small) indicates that the swimmer stays nearly straight.

In order to go further, and compute the (asymptotic) net displacement of the swimmer after one period of the oscillating external field, we linearize as well the equation \eqref{eq:position} to first order in $(\theta,\alpha_2,\alpha_3)$ near $(0,0,0)$. 

Noting, 
$$
\dot{x} = \vecG^x(\theta,\alpha_2,\alpha_3)  \cdot \begin{pmatrix} \dot{\theta} \\ \dot{\alpha}_2\\\dot{\alpha}_3\end{pmatrix}
\quad\textrm{and} \quad
\dot{y} = \vecG^y(\theta,\alpha_2,\alpha_3)  \cdot \begin{pmatrix} \dot{\theta} \\ \dot{\alpha}_2\\\dot{\alpha}_3\end{pmatrix}
$$
where $\vecG^x$ (resp. $\vecG^y$)  is the $1\times3$ matrix composed of $(\vecG_i\cdot\vece_x)_{i=1,\cdots3}$ (resp. $(\vecG_i\cdot\vece_y)_{i=1,\cdots3}$),
we obtain 
\begin{equation}
\small
\Delta x = \int_0^{\frac{2\pi}{\omega}}\epsilon  \left(\vecG^x(0,0,0)  +\epsilon\, { }^t\left(\tilde{\theta},\tilde{\alpha}_2,\tilde{\alpha}_3\right) \nabla\vecG^x(0,0,0)\right) \begin{pmatrix} \tilde{\theta}^\infty\\\tilde{\alpha}^\infty_2\\\tilde{\alpha}^\infty_3\end{pmatrix} \mbox{dt}' + o(\epsilon^2)\,.
\end{equation}

Since, $t\mapsto\begin{pmatrix} \tilde{\theta}^\infty(t)\\\tilde{\alpha}^\infty_2(t)\\\tilde{\alpha}^\infty_3(t)\end{pmatrix}$ is periodic, the latter equality reads
\begin{equation}
\label{eq:displacement}
\Delta x = \epsilon^2 \int_0^{\frac{2\pi}{\omega}} { }^t\left(\tilde{\theta}^\infty,\tilde{\alpha}^\infty_2,\tilde{\alpha}^\infty_3\right) \nabla\vecG^x(0,0,0) \begin{pmatrix} \tilde{\theta}^\infty\\\tilde{\alpha}^\infty_2\\\tilde{\alpha}^\infty_3\end{pmatrix} \mbox{dt}' + o(\epsilon^2)\,,
\end{equation}
and a straight forward computation leads to express 
\begin{equation*}
\small
\begin{aligned}
&\nabla\vecG^x(0,0,0) =\\&\frac{L}{2\left(2\eta+\eta_1\right)} \begin{pmatrix} 2(\eta-\eta_1)&-\eta_1 & \eta \\
-\frac{(6\eta\eta_1-4\eta\xi+\eta_1\xi_1}{(2\xi+\xi_1)}&-\frac{\eta_1(2\eta+\xi_1)}{(2\xi+\xi_1)}&-\frac{\eta(\eta_1-\xi_1)}{(2\xi+\xi_1)}\\
\frac{(2\eta^2+4\eta\eta_1-3\eta_1\xi)}{(2\xi+\xi_1)}&\frac{\eta_1(\eta-\xi)}{(2\xi+\xi_1)}&\frac{\eta(\eta+\eta_1+\xi)}{(2\xi+\xi_1)}\\
\end{pmatrix}\,.
\end{aligned}
\end{equation*}

Similarly, the same formula holds for $\Delta y$ by substituting $\vecG^x$ for $\vecG^y$ and in this case,
$$
\nabla\vecG^y(0,0,0) = \begin{pmatrix} 0 & 0& 0\\ 0  & 0 & 0\\0 & 0 & 0 \end{pmatrix}\,,
$$
thus, $\Delta y=o(\epsilon^2)$. It follows that the leading term, with respect to small angles, of the trajectory of the swimmer along the $y$-axis is negligible after one period of the oscillating fields compare to the one along the $x$-axis.

From now on, we focus on the $x$-displacement of the swimmer, $\Delta x$ and we prove that the leading term of order $\epsilon^2$ does not vanish after one period of the oscillating fields.

Plugging \eqref{eq:sol_edo_angles_infty} into \eqref{eq:displacement} and noting that the two terms vanish because of periodicity, 
\begin{equation*}
\small
\begin{aligned}
&\int_{0}^{\frac{2\pi}{\omega}} {}^t\left(\vecA^+(\omega) \exp{( i\omega t)} b\right) \nabla \vecG^x(0,0,0) \left( \vecA^+(\omega) \exp{( i\omega t)} b\right) = 0\,,\\
&\int_{0}^{\frac{2\pi}{\omega}} {}^t\left(\vecA^-(\omega) \exp{(- i\omega t)} b\right) \nabla \vecG^x(0,0,0) \left( \vecA^-(\omega) \exp{(- i\omega t)} b\right) = 0
\end{aligned}
\end{equation*} 
we obtain 
\begin{equation}
\small
\label{eq:delta_x_short}
\Delta x = \epsilon^2 \int_{0}^{\frac{2\pi}{\omega}} \frac{\omega}{4} \bigg[\left({}^t b {}^t\vecA^+(\omega)\left(\nabla \vecG^x(0,0,0) - {}^t\nabla\vecG^x(0,0,0)\right) \vecA^-(\omega) b\right)\bigg]\,.
\end{equation}

Notice already that since $A^+(0)=A^-(0)=-A^{-1}$, $\Delta x$ tends to $0$ when $\omega$ tends to $0$. A very low frequency produces no net motion (at order $\epsilon^2$), even after one period.

Moreover, the $3\times 3$ matrix $\left(\nabla \vecG^x(0,0,0 \,- {}^t\nabla\vecG^x(0,0,0)\right)$ is skew-symmetric and not null. Therefore, $0$ is an eigenvalue of multiplicity $1$. Let us denote by $\vecu$ its associated eigenvector. A direct computation, still using Mathematica, leads to 
$$
\vecu = 
\left(
\begin{array}{c}
\eta_1 \xi +\eta \xi_1-2 \eta \eta_1\\
2 \eta ^2+4 \eta_1 \eta -2 \xi  \eta -\xi_1 \eta -3 \eta_1 \xi \\
 6 \eta 
   \eta_1-2 \xi  \eta_1-4 \eta  \xi_1
\end{array}
\right)\,.
$$

Thus, to ensure that \eqref{eq:delta_x_short} is not null, it is sufficient to prove that the three of vectors $\left\{\vecu, \vecA^-(\omega) b, \vecA^+(\omega) b\right\}$ are independent.
But, for large frequencies $\omega$, we can expand the matrix $\vecA^{\pm}(\omega)$ as 
\begin{equation}
\label{eq:expansion_A}
\vecA^\pm(\omega) = \pm \frac{\vecI}{i\omega} -
\frac{\vecA}{\omega^2} +o(\frac{1}{\omega^2})\,.
\end{equation}
and 

\begin{equation*}
\small
\begin{aligned}
&\mbox{det}(\vecu, \vecA^-(\omega) b, \vecA^+(\omega) b) = \mbox{det}(\vecu,\frac{ \vecb}{\omega},\frac{\vecA \vecb}{\omega^2})  + o(\frac{1}{w^{4}})\\
 &=-\Bigg(L^9 \omega^3 \eta ^3 \eta_1^2 (8 \eta +7 \eta_1)^2 (3 \eta  \eta_1-2 \eta  \xi_1-\eta_1 \xi )\Bigg)^{-1} \times \\
 &\Bigg(\eta \eta_1 (113 K+38 M) + 216 M^2 (2 \eta +\eta_1) \bigl(\eta ^2 \xi_1 \bigl(4 \eta ^2 (5 K+2 M)\\
 & \quad\quad+\eta_1^2 (29 K+8 M)\bigr)+\eta \eta_1
   (\eta_1-\eta ) \bigl(K \left(4 \eta ^2+37 \eta \eta_1+13 \eta_1^2\right)\\\\
   &\quad\quad+6 \eta_1 M (2 \eta +\eta_1)\bigr)-\eta_1 \xi \bigl(8 \eta ^3 (2 K+M)\\
   &+2 \eta ^2 \eta_1(40 K+13 M)+\eta \eta_1^2 (53 K+14 M)+\eta_1^3 (13 K+6 M)\bigr)\bigr)\Bigg)\,.
\end{aligned}
\end{equation*}
This determinant does not vanish identically and we then obtain that by prescribing an oscillating magnetic field as \eqref{hyp:particular_field}, the magneto-elastic Purcell swimmer moves along the $x$-axis. Notice that here the assumption on the drag coefficients \eqref{drag_coeff_ass} is crucial.
\\

Moreover, $\Delta x$ tends to $0$ as $\omega \to 0$ and $\omega \to \infty$. This suggests the existence of an optimal frequency to drive the swimmer as was already observed in \cite{AlougesDeSimone14} (see Fig. $5$ and $7$).

\section{Discussion}
\label{4} 
This section underlines the challenge that we have to face in order to control this magnetic micro-swimmer. 
The previous result indicates that with a small sinusoidal magnetic field we are able to control the direction of the swimmer's displacement, but it does not imply neither global or local controllability properties. 

The latter local property is classically obtained by verifying the Kalman condition at an equilibrium point. 
Thus let us consider the system \eqref{eq:full_system}, around $(\vecX_e,\boldsymbol{u_e}) := \left((\vecx,0,0,0),(0,0)\right)$, which is an equilibrium point. At such a point, the swimmer is aligned with the horizontal field and thus, the torque due to the horizontal field vanishes leading to $\vecf_x(\vecX_e)=0$. Moreover, the system being invariant under translations, $\vecf_0$ does not depend on $(x,y)$, from which we deduce that the matrix range $R$ of the matrix $\nabla \vecf_0(\vecX_e)$ is at most of dimension 3. The Lie bracket 
$\mbox{ad}_{\vecf_0}(\vecf_x)(\vecX_e) = [\vecf_0,\vecf_x](\vecX_e)$ vanishes since $\vecf_0(\vecX_e)=0$ and $\vecf_x(\vecX_e)=0$.
By induction, for $k>1$
$$
\mbox{ad}^k_{\vecf_0}(\vecf_x)(\vecX_e) = [\vecf_0,\mbox{ad}^{k-1}_{\vecf_0}(\vecf_x)](\vecX_e) = 0\,.
$$
As far as $\mbox{ad}^k_{\vecf_0}(\vecf_y)(\vecX_e)$ is concerned, we have, still by induction, for all $k>1$
$$
\mbox{ad}^k_{\vecf_0}(\vecf_y)(\vecX_e) = [\vecf_0,\mbox{ad}^{k-1}_{\vecf_0}(\vecf_y)](\vecX_e) \in R\,.
$$
Therefore 
\begin{align*}
&\mbox{Span}\left(\left\{ \mbox{ad}_{\vecf_0}^k(\vecf_x)(\vecX_e),\mbox{ad}_{\vecf_0}^k(\vecf_y)(\vecX_e),\, k \geq 0\right\}\right) =\\& \mbox{Span}\left(\left\{\vecf_y(\vecX_e),\mbox{ad}_{\vecf_0}^k(\vecf_y)(\vecX_e),\, k \geq 0\right\}\right)\\
&\subset \mbox{Span}\left(\vecf_y(\vecX_e)\right)+R
\end{align*}
is at most of dimension 4. Thus the Kalman condition is not satisfied. It turns out that for non horizontal straight swimmers, i.e. $\vecX_e =  (\vecx,\theta,0,0)$ with $\theta\ne 0$, the same situation occurs due to the fact that 
\begin{equation}
\label{eq:f_x_null}
\vecf_x(\vecX_e)=-\tan(\theta)\vecf_y(\vecX_e).
\end{equation}

Moreover let us notice that
\begin{equation}
[\vecf_x,\vecf_y](\vecX_e)=f_{y,\vecX_e}^3 \vecf_y(\vecX_e)
\end{equation}
where $f_{y,\vecX_e}^3$ is the third component of the vector $\vecf_y(\vecX_e)$.\\
Therefore the Lie algebra $\text{\textbf{L}ie}\{\vecf_0,\vecf_x,\vecf_y\}|_{\vecX_e}$ span a vector space of dimension $4$. 
By changing the reference frame, a similar argument holds for all equilibrium points such as $\{(\vecx,\theta,0,0),\theta\in[0,2\pi]\}$.\\
This means that also the classical LARC condition is not satisfied and then the Sussmann condition does not hold (see \cite{Coron56}).\\

Here we have underlined the fact that the the magneto-elastic Purcell swimmer model is singular at the straight position which makes hard to get a local controllability result.
Of course if instead the swimmer starts at a non straight position, thanks to the boundedness of the drift, the magnetic field can compensate it and drive the swimmer.

\section{Conclusion}
In this paper,  
 by prescribing a particular oscillating field, we make an asymptotic expansion of the displacement of the swimmer, proving that this particular field allows to steer the swimmer along one direction.  Moreover we highlight the difficulties to get controllability result by showing that the classical conditions are not satisfied. It indicates that sophisticated techniques (see for instance \cite{GiraldiPomet15}) are required to obtain such controllability result.
%



\bibliography{N_Link_Magnetic}             
                                                   







\end{document}